\def\timestamp{%
Time-stamp: <non-image-of-Hstar.tex: Tuesday 25-09-2007 at 22:40:23 (cest)>}
\def\stripname Time-stamp: <#1 #2>{#2}
\edef\filedate{\expandafter\stripname\timestamp}

\documentclass{amsart}

\usepackage{amssymb,amsrefs}

\newcommand{\Ex}{\operatorname{Ex}}
\newcommand{\Int}{\operatorname{int}}
\newcommand{\orpr}[2]{\langle#1,#2\rangle}
\newcommand{\preim}{^\gets}

\DeclareMathSymbol\diamond0{AMSa}{"06}
\DeclareMathSymbol\restr2{AMSa}{"16}
\DeclareMathSymbol\A0{AMSb}{`A}
\DeclareMathSymbol\Dow0{AMSb}{`D}
\DeclareMathSymbol\HH0{AMSb}{`H} 
\DeclareMathSymbol\LL0{AMSb}{`L} 
\DeclareMathSymbol\N0{AMSb}{`N} 
\DeclareMathSymbol\R0{AMSb}{`R} 
\DeclareMathSymbol\mapdiagbin{2}{AMSa}{"4D}

\newcommand\betaH{\beta\HH}
\newcommand\betaN{\beta\N}
\newcommand\Hstar{\HH^*}
\newcommand\Nstar{\N^*}

\let\phi\varphi

\newcommand{\CH}{\mathsf{CH}}
\newcommand{\MA}{\mathsf{MA}}
\newcommand{\OCA}{\mathsf{OCA}}

\newcommand{\cee}{\mathfrak{c}}

\newcommand{\Pow}{\mathcal{P}}
\newcommand{\fin}{\mathit{fin}}

\newtheorem{theorem}{Theorem}[section]
\newtheorem{lemma}[theorem]{Lemma}    
\newtheorem{proposition}[theorem]{Proposition}    

\theoremstyle{remark}

\newtheorem{question}[theorem]{Question}

\begin{document}

\title{A separable non-remainder of $\HH$}

\author[Alan Dow]{Alan Dow\dag}
\address{
Department of Mathematics\\
UNC-Charlotte\\
9201 University City Blvd. \\
Charlotte, NC 28223-0001}
\email{adow@uncc.edu}
\urladdr{http://www.math.uncc.edu/\~{}adow}
\thanks{\dag Supported by NSF grant DMS-0554896}

\author{Klaas Pieter Hart}
\address{Faculty of Electrical Engineering, Mathematics and Computer Science\\
         TU Delft\\
         Postbus 5031\\
         2600~GA {} Delft\\
         the Netherlands}
\email{k.p.hart@tudelft.nl}
\urladdr{http://fa.its.tudelft.nl/\~{}hart}

\date{\filedate}

\begin{abstract}
We prove that there is a compact separable continuum that
(consistently) is not a remainder of the real line.
\end{abstract}

\subjclass[2000]{Primary: 54F15. 
                 Secondary: 03E50, 03E65, 54A35, 54D15, 54D40, 54D65}
\keywords{separable continuum, continuous image, $\Hstar$, $\beta X$,  
          $\OCA$}

\maketitle

\section*{Introduction}

Much is known about the continuous images of $\Nstar$, the \v{C}ech-Stone
remainder of the discrete space~$\N$.
It is nigh on trivial to prove that every separable compact Hausdorff space
is a continuous image of~$\Nstar$ (we abbreviate this as `$\Nstar$-image'), 
it is a major result of Parovi\v{c}enko, from~\cite{Parovicenko63},
that every compact Hausdorff space of weight~$\aleph_1$ is an $\Nstar$-image
and in~\cite{MR671232} Przymusi\'nski 
used the latter result to prove that all perfectly normal
compact spaces are $\Nstar$-images.
Under the assumption of the Continuum Hypothesis Parovi\v{c}enko's result 
encompasses all three results: a compact Hausdorff space is an $\Nstar$-image
if and only if it has weight~$\cee$ or less.

In~\cite{MR1707489} the authors formulated and proved a version of 
Parovi\v{c}enko's theorem in the class of continua:
every continuum of weight $\aleph_1$ is a continuous image of~$\Hstar$
(an `$\Hstar$-image'),
the \v{C}ech-Stone remainder of the subspace~$\HH=[0,\infty)$ of the real line.
This result built on and extended the corresponding result for metric continua
from~\cite{MR0214033}.
Thus the Continuum Hypothesis ($\CH$) allows one to characterize the 
$\Hstar$-images as the continua of weight~$\cee$ or less.
The paper~\cite{MR1707489} contains further results on $\Hstar$-images that 
parallel older results about $\Nstar$-images:
Martin's Axiom ($\MA$) implies all continua of weight less than~$\cee$ 
are $\Hstar$-images, in the Cohen model the long segment of length~$\omega_2$
is not an $\Hstar$-image, and it is consistent with $\MA$ that not every 
continuum of weight~$\cee$ is an $\Hstar$-image.

The natural question whether the `trivial' result on separable compact spaces
has its parallel version for continua proved harder to answer than expected.
We show that in this case the parallelism actually breaks down.
There is a well-defined separable continuum~$K$ that is not an $\Hstar$-image
if the Open Colouring Axiom ($\OCA$) is assumed.
This also answers a more general question raised by G.~D. Faulkner
(\cite{MR1707489}*{Question~7.3}):
if a continuum is an $\Nstar$-image must it be an $\Hstar$-image?
Indeed, $K$~\emph{is} separable and hence an $\Nstar$-image.

It is readily seen that $\betaH$ itself is an $\Hstar$-image: by moving 
back and forth in ever larger sweeps one constructs a map from~$\HH$ onto 
itself whose \v{C}ech-Stone extension maps~$\Hstar$ onto~$\betaH$.
Indeed the same argument applies to any space that is the union of
a connected collection of Peano continua: its \v{C}ech-Stone compactification
is an~$\Hstar$-image.
Thus, e.g., for every~$n$ the space~$\beta\R^n$ is an $\Hstar$-image.
Our example is one step up from these examples: it is the \v{C}ech-Stone
compactification of a string of $\sin\frac1x$-curves.

Our result also shows that the proof in~\cite{MR1707489} cannot be extended
beyond~$\aleph_1$, as $\OCA$ is compossible with Martin's Axiom ($\MA$).
The adage that $\MA$ makes all cardinals below~$\cee$ behave as if they
are countable would suggest that the aforementioned proof, an inverse-limit
construction, could be made~$\cee$ long, at least if $\MA$ holds.
We see that this is not possible, even if the continuum is separable.

The paper is organized as follows.
Section~\ref{sec.prelim} contains a few preliminaries, including the 
consequences of $\OCA$ that we shall use.
In Section~\ref{sec.thespace} we construct the continuum~$K$ and show how
$\OCA$~implies that it is not an $\Hstar$-image.
Finally, in Section~\ref{sec.remarks} we give a few more details on the
lack of efficacy of~$\MA$ in this and we discuss and ask whether other 
potential $\Hstar$-images are indeed $\Hstar$-images.

We thank the referee for pointing out that there was much room for improvement
in our presentation.

\section{Preliminaries}
\label{sec.prelim}

\subsection*{Closed and open sets in $\beta X$}

Since we will be working with subsets of the plane we can economize
a bit on notation and write $\beta F$ for the closure-in-$\beta X$
of a closed subset of the space~$X$ itself; 
we also write $F^*=\beta F\setminus F$ .
If $O$~is an open subset of~$X$ then 
$\Ex O=\beta X\setminus\beta(X\setminus O)$
is the largest open subset of~$\beta X$ whose intersection with~$X$ is~$O$.

In dealing with closed subsets of $\Hstar$ the following, 
which is Proposition~3.2 from~\cite{Hart}, is very useful.

\begin{proposition}\label{prop.nice.open}
Let $F$ and $G$ be disjoint closed sets in\/~$\Hstar$.
There is an increasing and cofinal sequence sequence 
$\langle a_k:k\in\omega\rangle$ in~$\HH$ such that 
$F\subseteq\Ex\bigcup_k(a_{2k+1},a_{2k+2})$ and 
$G\subseteq\Ex\bigcup_k(a_{2k},a_{2k+1})$. \qed
\end{proposition}

We shall be working with closed subsets of the plane (or~$\HH$) that can be 
written as the union of a discrete sequence $\langle F_n:n\in\omega\rangle$ 
of compact sets.
The extension of the natural map~$\pi$ from~$F=\bigcup_nF_n$ to~$\omega$,
that sends the points of~$F_n$ to~$n$, partitions $\beta F$ into sets
indexed by~$\beta\omega$: for $u\in\beta\omega$ we write 
$F_u=\beta\pi\preim(u)$.
If the $F_n$ are all connected then so is every~$F_u$
and, indeed, the $F_u$ are the components of~$\beta F$, 
see \cite{Hart}*{Corollary 2.2}.

For use below we note the following.

\begin{lemma}\label{lemma.irr}
If each $F_n$ is an irreducible continuum, between the points~$a_n$ and~$b_n$
say, then so is each~$F_u$, between the points~$a_u$ and~$b_u$. \qed
\end{lemma}

\subsection*{The Open Colouring Axiom}

The \emph{Open Coloring Axiom} ($\OCA$) was formulated by
To\-dor\-\v{c}e\-vi\'c in~\cite{MR980949}.
It reads as follows:
if $X$ is separable and metrizable and if $[X]^2=K_0\cup K_1$,
where $K_0$~is open in the product topology of~$[X]^2$, then
\emph{either} $X$~has an uncountable $K_0$-homogeneous subset~$Y$
\emph{or} $X$~is the union of countably many $K_1$-homogeneous subsets.

One can deduce the conjunction $\OCA$ and $\MA$ from the 
\emph{Proper Forcing Axiom} or prove it consistent in an $\omega_2$-length 
countable support proper iterated forcing construction, 
using~$\diamond$ on~$\omega_2$ to predict all possible subsets of the Hilbert 
cube and all possible open colourings of these, as well as all possible ccc
posets of cardinality~$\aleph_1$.

We shall make use of $\OCA$ only but we noted the compossibility with $\MA$
in order to substantiate the claim that the latter principle does not imply
that all separable continua are $\Hstar$-images.

\subsubsection*{Triviality of maps}

We shall use two consequences of $\OCA$.
The first says that continuous surjections from $\omega^*$ onto~$\beta\omega$
are `trivial' on large pieces of~$\omega^*$.
If $\phi:\omega^*\to\beta\omega$ is a continuous surjection then it induces,
by Stone duality, an embedding of $\Phi:\Pow(\omega)\to\Pow(\omega)/\fin$
by $\Phi(A)=\phi\preim[\beta A]$.
The following is a consequence of~\cite{MR1752102}*{Theorem~3.1}, where
for a subset $M$ of~$\omega$ we write $\tilde M=(M-1)\cup M\cup(M+1)$.

\begin{proposition}[$\OCA$]\label{prop.triviality}
With the notation as above there are 
infinite subset $D$ and $M$ of~$\omega$ and a map $\psi:D\to\tilde M$
such that for every subset $A$ of $\tilde M$ one has 
$\Phi(A)=\psi\preim[A]^*$. \qed
\end{proposition}

Thus, on the set $D^*$ the map $\phi$ is determined by the map 
$\psi:D\to\tilde M$; this is the sense in which $\phi\restr D^*$ might
be called trivial.
It is also important to note that $D=^*\psi\preim[\tilde M]$,
which follows from $D^*=\phi\preim[\beta\tilde M]$;
this will be used in our proof.

\subsubsection*{Non-images of $\Nstar$}
 
The final nail in the coffin of a purported map from $\Hstar$ onto the 
continuum~$K$ will be the following result from~\cite{MR1679586},
where $\Dow=\omega\times(\omega+1)$.

\begin{proposition}[$\OCA$]\label{prop.Dow.not.image}
The \v{C}ech-Stone remainder\/ $\Dow^*$ is not an $\Nstar$-image. \qed 
\end{proposition}

\section{The non-image}
\label{sec.thespace}

\subsection*{The example}

We start by replicating the $\sin\frac1x$-curve along the $x$-axis in the 
plane: for $n\in\omega$ we set 
$K_n=\bigl(\{n\}\times[-1,1]\bigr)\cup
     \bigl(\bigl\{\orpr{n+t}{\sin\frac\pi{t}}:0< t\le1\bigr\}\bigr)$.
The union $K=\bigcup_nK_n$ is connected and its \v{C}ech-Stone
compactification $\beta K$ is separable continuum.
We shall show that $\OCA$ implies that $\beta K$ is not a continuous image
of~$\Hstar$.

We define four closed sets that play an important part in the proof.
For $n\in\omega$ we define:
\begin{itemize}
\item $S_n=\{\orpr xy\in X:n+\frac13\le x\le n+\frac23\}$,
\item $S_n^+=\{\orpr xy\in X:n+\frac14\le x\le n+\frac34\}$,
\item $T_n=\{\orpr xy\in X:n-\frac14\le x\le n+\frac14\}$,
\item $T_n^+=\{\orpr xy\in X:n-\frac13\le x\le n+\frac13\}$;
\end{itemize}
we put $S=\bigcup_{n\in\omega}S_n$, $S^+=\bigcup_{n\in\omega}S_n^+$,
$T=\bigcup_{n\in\omega}T_n$ and $T^+=\bigcup_{n\in\omega}T_n^+$.
Note that $S\cap T=\emptyset$ and hence $\beta S\cap\beta T=\emptyset$ 
in $\beta X$.

We also note that the four sets $S_n$, $S_n^+$, $T_n$ and $T_n^+$ are
all connected and that we therefore know exactly what the components
of $\beta S$, $\beta S^+$, $\beta T$ and $\beta T^+$ are.
Note that, by Lemma~\ref{lemma.irr}, each continuum $T_u^+$ (as well as
$S_u$, $S_u^+$ and $T_u$) is irreducible, as each $T_n^+$ is irreducible
(between its end~points $\orpr{n-\frac13}{-1}$ and $\orpr{n+\frac13}{0}$).

Finally we  note that $S_n$ meets the sets $T_n^+$ and $T_{n+1}^+$ only,
that $T_n$ meets $S_{n-1}^+$ and $S_n^+$ only, etcetera.
This behaviour persists when we move to the continua~$S_u$, $T_u^+$, 
$S_u^+$ and~$T_u^+$, when we define $u+1$ and $u-1$, for $u\in\omega^*$,
in the obvious way: $u+1$~is generated by~$\{A+1:A\in u\}$ and $u-1$~is
generated by $\{A:A+1\in u\}$.

\subsection*{Properties of a potential surjection}

Assume $h:\Hstar\to\beta K$ is a continuous surjection and apply 
Proposition~\ref{prop.nice.open} to the closed subsets
$h\preim[\beta S]$ and $h\preim[\beta T]$ of~$\Hstar$
to get a sequence $\langle a_k:k\in\omega\rangle$.
After composing $h$ with a piecewise linear map we may assume, without
loss of generality, that $a_k=k$ for all~$k$.
We obtain  
$h\preim[\beta S]\cap\beta\bigcup_{k\in\omega}I_{2k+1}=\emptyset$ and
$h\preim[\beta T]\cap\beta\bigcup_{k\in\omega}I_{2k}=\emptyset$, where
$I_k=[k,k+1]$.
We write $2\omega$ and $2\omega+1$ for the sets of even and odd natural
numbers respectively.

The map $h$ induces maps from $(2\omega)^*$ and $(2\omega+1)^*$ 
onto $\beta\omega$, as follows.
If $u\in(2\omega)^*$ then $h[I_u]$ is a connected set that is disjoint 
from $\beta T$, hence it must be contained in a component of $\beta S^+$.
Likewise, if $v\in(2\omega+1)^*$ then $h[I_v]$ is contained in a component
of $\beta T^+$.
Thus we get maps $\phi_0:(2\omega)^*\to\beta\omega$ and 
$\phi_1:(2\omega+1)^*\to\beta\omega$ defined by
\begin{itemize}
\item $\phi_0(u)=x$ iff $h[I_u]\subseteq S_x^+$,
\item $\phi_1(v)=y$ iff $h[I_v]\subseteq T_y^+$.
\end{itemize}

\begin{lemma}
The maps $\phi_0$ and $\phi_1$ are continuous.  
\end{lemma}

\begin{proof}
For $k\in\omega$ put $r_k=k+\frac12$.
Observe that, by connectivity,  $h(r_u)\in S_x^+$ iff 
$h[I_u]\subseteq S_x^+$, so that $\phi_0$ can be decomposed as
$u\mapsto r_u\mapsto h(r_u)\mapsto \pi_0(h(r_u))$, where
$\pi_0:\beta S^+\to\beta\omega$ is the natural map.

The argument for $\phi_1$ is similar. 
\end{proof}

The maps $\phi_0$ and $\phi_1$ are not unrelated.
Let $u\in(2\omega)^*$ and put $x=\phi_0(u)$.
Then $h[I_u]\subseteq S_x^+$, so that $h[I_{u-1}]$ and $h[I_{u+1}]$
both intersect $S_x^+$.
However, $S_x^+$ intersects only the continua $T_x^+$ and $T_{x+1}^+$,
so that $\phi_1(u-1),\phi_1(u+1)\in\{x,x+1\}$.
By symmetry a similar statement can be made if $y=\phi_1(v)$: then
$\phi_0(v-1),\phi_0(v+1)\in\{y-1,y\}$.

Using these relationships we can deduce some extra properties of $\phi_0$
and $\phi_1$.

\begin{lemma}\label{lemma.phi0}
If $u\in(2\omega)^*$ then 
$\phi_0(u-2)$ and $\phi_0(u+2)$ both are in $\{x-1,x,x+1\}$, 
where $x=\phi_0(x)$. \qed
\end{lemma}

\section{An application of $\OCA$}

We apply Proposition~\ref{prop.triviality} to the embedding $\Phi_0$ of
$\Pow(\omega)$ into $\Pow(2\omega)/\fin$ defined by 
$\Phi_0[A]=\phi_0\preim[\beta A]$.
We find infinite sets $D\subseteq2\omega$ and $M\subseteq\omega$
together with a map 
$\psi:D\to \tilde M$ that induces~$\Phi_0$ on its range:
for every subset $A$ of $\tilde M$ we have
$\Phi_0[A]=\psi_0\preim[A]^*$.
As noted above this implies that $\phi_0\restr D^*=\beta\psi_0\restr D^*$.

For $m\in\tilde M$ and $u\in(2\omega)^*$ we have the equivalence
$\phi_0(u)=m$ iff $\psi_0\preim[\{m\}]\in u$.
Using the properties of $\phi_0$ stated in Lemma~\ref{lemma.phi0}
we deduce the following inclusion-mod-finite
$$
\bigl(\psi_0\preim\bigl[\{m\}\bigr]-2\bigr)\cup 
      \psi_0\preim\bigl[\{m\}\bigr]\cup
\bigl(\psi_0\preim\bigl[\{m\}\bigr]+2\bigr)
\subseteq^*
\psi_0\preim\bigl[\{m-1,m,m+1\}\bigr]
\subseteq D
$$ 
Therefore we get for every $m\in M$ a $j_m$ such that:
if $n\ge j_m$ and $\psi_0(n)=m$ then $n-2,n+2\in D$.

\begin{lemma}
For every $m\in M$ there are infinitely many $n\in D$ such that $\psi_0(n)=m$ 
and $\psi_0(n+2)\neq m$.  
\end{lemma}

\begin{proof}
Let $m\in M$ and take $m'\in M\setminus\{m\}$.
Let $n\in D$ be arbitrary such that $\psi_0(n)=m$ and $n\ge j_m$;
choose $n'>n$ such that $\psi_0(n')=m'$.
There must be a first index $i$ such that $\psi_0(n+2i)\neq\psi_0(n+2i+2)$ as
otherwise we could show inductively that $n+2i\in D$ and $\psi_0(n+2i)=m$
for all~$i$, which would imply that $\psi_0(n')=m$.
For this minimal $i$ we have $\psi_0(n+2i)=m$ and $\psi_0(n+2i+2)\neq m$.
\end{proof}

We use this lemma to find an infinite subset~$L$ of~$D$ where $\phi_0$
and $\phi_1$ are very well-behaved. 

Let $m_0=\min M$ and choose $l_0\ge j_{m_0}$ such that $\psi_0(l_0)=m_0$
and $\psi_0(l_0+2)\neq m_0$.
Proceed recursively: choose $m_{i+1}\in M$ larger than $m_i+3$
and $\psi_0(l_i+2)+3$, and then pick $l_{i+1}$ larger
than~$l_i$ and $j_{m_{i+1}}$ such that $\psi_0(l_{i+1})=m_{i+1}$
and $\psi_0(l_{i+1}+2)\neq m_{i+1}$.

Consider the set $L=\{l_i:i\in\omega\}$ and thin out $M$ so that it will
be equal to $\{m_i:i\in\omega\}$.
Let $u\in L^*$ and let $x=\phi_0(u)=\psi_0(u)$; we assume, 
without loss of generality, that $\{l\in L:\psi_0(l+2)=\psi_0(l)+1\}$
belongs to~$u$.
It follows that $\phi_0(u+2)=x+1$ and this means that $\phi_1(u+1)=x$.

We find that
$h[I_u]\subseteq S_x^+$,
$h[I_{u+1}]\subseteq T_x^+$ and $h[I_{u+2}]\subseteq S_{x+1}^+$.
Therefore the image $h[I_u\cup I_{u+1}\cup I_{u+2}]$ is a subcontinuum of
$S_x^+\cup T_x^+\cup S_{x+1}^+$ that meets $S_x^+$ and $S_{x+1}^+$.
Because $T_x^+$ is irreducible we find that 
$T_x^+\subseteq h[I_u\cup I_{u+1}\cup I_{u+2}]$ and hence
$T_x\subseteq h[I_{u+1}]$,
because the other two parts of this continuum are disjoint from~$\beta T$.

We now have infinite sets $L$ and $M$ where the maps $\phi_0$ and $\phi_1$
behave very nicely indeed.
Because $\psi_0$ maps $L$ onto~$M$ the map $\phi_0$ maps $L^*$ onto~$M^*$.
Furthermore, if $u\in L^*$ and $x=\phi_0(u)$ then also $x=\phi_1(u+1)$
and $T_x\subseteq h[I_{u+1}]\subseteq T_x^+$.

We put $\LL=\bigcup\{I_{u+1}:u\in L^*\}$ and we observe that, by the
inclusions above, 
$$
\bigcup\{T_x:x\in M^*\}\subseteq h[\LL]\subseteq \bigcup\{T_x^+:x\in M^*\}.
\eqno(*)
$$
We put $h_L=h\restr \LL$.

\subsection*{A map from $\Nstar$ onto $\Dow^*$}

We now use $h_L$ to create a map from $\Nstar$ onto $\Dow^*$, which will 
yield the contradiction that finishes the proof.

Let $F=\bigcup_{m\in M}T_m\cap\bigl([0,\infty)\times[\frac12,1]\bigr)$ and 
$G=\bigcup_{m\in M}T_m\cap\bigl([0,\infty)\times[0,1]\bigr)$.
Observe that the inclusion map from $F$ to $G$ induces the identity map
between their respective component spaces and hence also the identity
map between the component spaces of $\beta F$ and $\beta G$.
We work with the closed subsets~$F^*$ and~$G^*$ of~$K^*$; the former is
contained in the interior of the latter, hence the same holds for
$h_L\preim[F^*]$ and $h_L\preim[G^*]$.
We apply Proposition~\ref{prop.nice.open} and obtain, for every $l\in L$, 
a finite family $\mathcal{I}_l$ of subintervals of~$I_l$ such that for the 
closed set $H=\bigcup_{l\in L}\bigcup\mathcal{I}_l$ we have
$$
h_L\preim[F^*]\subseteq \Int H^*\subseteq H^*\subseteq\Int h_L\preim[G^*]
\eqno(\dag)
$$
Endow the countable set of intervals $\mathcal{I}=\bigcup_{l\in L}\mathcal{I}_l$
with the discrete topology and let $p\in\mathcal{I}^*$; the corresponding
component of~$H^*$ is mapped by~$h_L$ into a component of~$G^*$.
Thus we obtain a map from $\mathcal{I}^*$ into the component space of $G^*$.
This map is onto: let $C_G$ be a component of $G^*$ and let $C_F$ be the
unique component of~$F^*$ contained in~$C_G$.
Because of~$(\dag)$ and $(*)$ there is a family of components of~$H^*$
that covers~$C_F$; all these components are mapped into~$C_G$.

We obtain a map from $\mathcal{I}^*$ onto the component space of $G^*$.
This map is continuous; this can be shown as for the maps~$\phi_0$ and~$\phi_1$
using midpoints of the intervals and the quotient map from $G^*$ onto its 
component space.
The component space of~$G$ itself is $\Dow$, so that $G^*$ has $\Dow^*$ as its
component space.
Thus the assumption that $\Hstar$ maps onto $\beta K$ leads, assuming $\OCA$,
to a continuous surjection from $\omega^*$ onto~$\Dow^*$, which, 
by Proposition~\ref{prop.Dow.not.image} is impossible.

\section{Further remarks}
\label{sec.remarks}

\subsection{Comments on the construction}

The proofs in \cites{MR1679586,MR1752102} that certain spaces are \emph{not}
$\Nstar$-images follow the same two-step pattern: first show that no `trivial'
map exists and then show that $\OCA$ implies that if there is a map at all
then there must also be a `trivial' one.
In the context of our example it should be clear that there is no map
from~$\HH$ to the plane that induces a map from~$\Hstar$ onto~$\beta K$;
it would have been nice to have found a map from $\bigcup_{l\in L}I_{l+1}$
to the plane that would have induced~$h_L$ but we did not see how to construct
one.

\subsection{$\MA$ is not strong enough}

As mentioned in the introduction the principal result of~\cite{MR1707489} 
states that every continuum of weight~$\aleph_1$ is an $\Hstar$-image.
In that paper the authors also prove that under $\MA$ every continuum of 
weight less than~$\cee$ is an $\Hstar$-image; the starting point of that
proof was the result of Van Douwen and Przymusi\'nski \cite{MR561588}
that, under $\MA$, every compact Hausdorff space of weight less than~$\cee$
is an~$\Nstar$-image.
Given such a continuum~$X$, of weight~$\kappa<\cee$, one assumes it is 
embedded in the Tychonoff cube~$I^\kappa$ and takes a continuous 
map $f:\betaN\to I^\kappa$ such that $f[\Nstar]=X$.
What the proof then establishes, using~$\MA$, is that $f$ has an extension
$F:\betaH\to I^\kappa$ such that $F[\Hstar]=X$.
Thus, in a very real sense, one can simply connect the dots of~$\N$ to
produce a map from~$\Hstar$ onto~$X$ that extends the given map from~$\Nstar$
onto~$X$.

Since $\MA$ and $\OCA$ are compossible our example shows that $\MA$ does not
imply that all separable continua are $\Hstar$-images and, a fortiori, 
that the two proofs from~\cite{MR1707489} cannot be amalgamated to show that 
the answer to Faulkner's question is positive under~$\MA$, not even for 
separable spaces.

\subsection{Other images}

As noted in the introduction there are many parallels between
the results on $\Nstar$-images and those on $\Hstar$-images.
The example in this paper shows that there is no complete parallelism.
There are some results on $\Nstar$-images where no parallel has been found
or disproved to exist.

We mentioned Przymusi\'nski's theorem from~\cite{MR671232} that every 
perfectly normal compact space is an $\Nstar$-image.
By compactness every perfectly normal compact space is first-countable
and by Arhangel\cprime ski\u{\i}'s theorem (\cite{MR0251695}) 
every first-countable compact space has weight~$\cee$ and is therefore
an $\Nstar$-image if $\CH$ is assumed.

Thus we get two questions on $\Hstar$-images.

\begin{question}
Is every perfectly normal compact continuum an $\Hstar$-image?
\end{question}

\begin{question}
Is every first-countable continuum an $\Hstar$-image?  
\end{question}

The questions are related of course but the question on first-countable
continua might get a consistent negative answer sooner than the one on
perfectly normal continua in the light of Bell's consistent example, 
from~\cite{MR1058795},
of a first-countable compact space that is not an~$\Nstar$-image.

\end{document}